\documentclass[letterpaper, 10 pt, conference]{ieeeconf}  

\IEEEoverridecommandlockouts                              

\usepackage{amsfonts,bm}
\usepackage{amssymb,amsmath}
\usepackage{graphicx}
\usepackage{cite}
\usepackage{algorithm2e}
 \RestyleAlgo{ruled}

\renewcommand{\epsilon}{\varepsilon}
\renewcommand{\phi}{\varphi}

\renewcommand{\d}{\,\mathrm{d}}

\newtheorem{proposition}{Proposition}
\newtheorem{theorem}{Theorem}

\newtheorem{corollary}{Corollary}
\newtheorem{lemma}{Lemma}

\renewcommand{\epsilon}{\varepsilon}
\renewcommand{\phi}{\varphi}

\newcommand{\R}{\mathbb{R}}

\title{\LARGE \bf
Exact Cost‑Increment Formula for Optimal Control of Semilinear Evolution Equations$^*$
}

\author{Roman Chertovskih, Nikolay Pogodaev, Maxim Staritsyn, and A. Pedro Aguiar, \IEEEmembership{Member, IEEE}-
\thanks{$^*$RC, MS and PA acknowledge financial support of the 
Research Center for Systems and Technologies 
(UID/00147), the Associate Laboratory ARISE – 
(LA/P/0112/2020, DOI: 10.54499/LA/P/0112/2020)
and the project 2023.09597.CBM, all funded by Fundação para a Ciência e a Tecnologia, I.P./MECI through national funds. 
NP is supported by a subsidy from the Ministry of Education and Science of Russia (project no.~121041300060-4).
}
\thanks{Roman Chertovskih, Maxim Staritsyn, and A. Pedro Aguiar are with Research Center for Systems and Technologies (SYSTEC), ARISE, Department of Electrical
and Computer Engineering,
Faculdade de Engenharia, Universidade do Porto, Rua Dr. Roberto Frias, s/n 4200-465, Porto, Portugal (e-mail: roman@fe.up.pt, staritsyn@fe.up.pt, pedro.aguiar@fe.up.pt). Matrosov Institute for System Dynamics and
Control Theory (ISDCT~SB~RAS), Irkutsk, Russia (e-mail: nickpogo@gmail.com)
}%
}

\begin{document}
\maketitle
\thispagestyle{empty}
\pagestyle{empty}

\begin{abstract}
We address optimal control of semilinear evolution equations on Banach spaces with finitely many control channels, a framework encompassing a broad class of infinite‑dimensional dynamical systems, arising in many applications. For this setting, we derive an exact and global formula quantifying the increment of the cost functional with respect to an arbitrary reference control. This identity enables the design of monotone descent algorithms that require no linearization or step‑size tuning. We further establish the existence of optimal controls and propose a practical sample‑and‑hold realization of the descent step suitable for numerical implementation. The effectiveness of the method is demonstrated on a controlled reaction–diffusion equation.\end{abstract}

\section{INTRODUCTION}\label{sec:I}

Semilinear evolution equations on Banach spaces provide a common framework for many distributed control systems, including those governed by parabolic and hyperbolic PDEs, delay and integro-differential equations, or coupled systems \cite{bensoussan2007representation,Fattorini1999,Li1995,trelatControlFiniteInfinite2024}. In this setting, first-order optimality conditions are classical, while numerical optimization methods are often based on local linearization and step-size selection.

This paper takes a different route. Fixing a reference control, we derive an identity that expresses the full increment of the cost functional for an arbitrary admissible control through its deviation from the reference one --- an \emph{exact cost-increment formula}. Such identities are useful in two rather different ways. On the one hand, they fit into the line of nonclassical necessary conditions developed in \cite{Dykht2022,chertovskihOptimalControlDiffusion2024,SChPP-2022}. On the other hand, they lead naturally to monotone improvement procedures, close in spirit to the numerical methods of \cite{srochko1982computational,borzi2023sequential}. In the present paper, we pursue only this second direction.

The main contribution of this work is to generalize the variational framework developed in \cite{chertovskihOptimalControlDiffusion2024,SChPP-2022,pogodaevExactFormulaeIncrement2024,CPSA2023} to the setting of semilinear evolution equations on Banach spaces covering a broad class of infinite‑dimensional control systems. Section~II introduces the control problem and the standing assumptions, establishes the regularity properties of the mild flow needed for the backward cost construction, and proves existence of minimizers under a finite-channel assumption on the control operator. Section~III derives an exact formula for the full cost increment relative to a baseline control, shows how this formula generates a monotone descent mechanism, and discusses a simple sample-and-hold realization of the descent step. The method is illustrated there on a reaction--diffusion example. The proofs of the main results are collected in the appendix.

\section{Problem Statement and Analytical Foundations}

Let $X \doteq (X,\|\cdot\|_{X})$ be a real Banach space, and $T>0$, $\alpha\geq 0$ be given. On a fixed time interval $I \doteq [0,T]$, consider the optimal control problem
\begin{equation*}
({\bf P}) \quad
\inf\Big\{\mathcal I[u] \doteq \ell(x_T^u) + \frac{\alpha}{2}\int_I\bigl|u(t)\bigr|^{2} \d t\colon \ u \in \mathcal U\Big\}.
\end{equation*}
Here $\ell\colon X \to \mathbb R$ is a terminal cost functional. The admissible controls form the set $\mathcal U \doteq L^\infty(I;U)$, viewed as a subset of $L^\infty(I;\R^m)$, where $U \subset \R^m$ is compact and convex. For each $u\in\mathcal U$, the corresponding state trajectory $x \doteq x^u\colon I \to X$, $t \mapsto x_t$, is governed by the semilinear evolution equation: 
\begin{equation}
\dot x_t = A x_t + F_t(x_t,u(t)) \doteq A x_t + f_t(x_t) + G_t(x_t)\,u(t) \label{ODE}
\end{equation}
for almost each (a.e.) $t \in I$, with prescribed initial condition $x_0 \in X$. The operator $A\colon D(A)\subset X \to X$ is assumed to generate a $C_0$ semigroup on $X$, whereas the maps $f\colon I \times X \to X$ and $G\colon I \times X \to \mathcal L(\R^m;X)$ are given. 

Our choice of the class $\mathcal U$ puts $(\mathbf P)$ in the realm of \emph{ensemble control}, i.e., a centralized control of distributed systems.

\subsection{Notation, Standing Assumptions \& Preliminaries}\label{ssec:assumptio}

Let $Y \doteq (Y,\|\cdot\|_{Y})$ be a Banach space. Measurability of a map $\phi\colon I \to Y$ is always understood in the strong (Bochner) sense. 
By $C(I;Y)$ we denote the space of continuous maps $\phi\colon I \to Y$ with the norm $\|\phi\|_\infty \doteq \sup_{t \in I}\|\phi_t\|_{Y}$.

If $Z$ is another Banach space, then $\mathcal L(Y;Z)$ stands for the space of bounded linear maps $L\colon Y \to Z$, equipped with the operator norm $\|L\|_{\mathcal L(Y;Z)} \doteq \sup_{\|\mathrm y\|_{Y} \leq 1}\|L\mathrm y\|_{Z}$. The adjoint of $L$ is denoted by $L'\colon Z' \to Y'$. By $C^1(Y;Z)$ we mean the class of continuously Fr\'echet differentiable maps $F\colon Y \to Z$, that is, such that $DF\colon Y \to \mathcal L(Y;Z)$ is continuous. If $Z=\R$, we simply write $C^1(Y)$. 
The identity map on a space $Y$ is denoted by $\mathrm{id}_Y$.

If $A\colon D(A)\subset Y \to Y$ is a linear operator, then $D(A)$ denotes its domain. A strongly continuous semigroup on $Y$ is a family $(S_t)_{t\geq 0}\subset \mathcal L(Y;Y)$ such that $S_0=\mathrm{id}_Y$, $S_{t+\tau}=S_t\circ S_\tau$ for all $t,\tau\geq 0$, and the orbit map $t\mapsto S_t\mathrm x$ is continuous for every $\mathrm x\in Y$.

We impose the following assumptions $(\mathbf{A})$:\footnote{We intentionally formulate all hypotheses with a margin. These are not the weakest conditions under which the arguments work, but they keep the proofs more transparent.}
For every \(\mathrm{x}\in X\), the maps \(t\mapsto f_t(\mathrm{x})\) and \(t\mapsto G_t(\mathrm{x})\) are measurable. There exist constants \(M_f,M_G\geq 0\) such that, for all \(\mathrm{x},\mathrm{y}\in X\) and for a.a. \(t\in I\), one has
\[
\begin{gathered}
\|f_t(\mathrm{x})-f_t(\mathrm{y})\|_X\le M_f\,\|\mathrm{x}-\mathrm{y}\|_X,\\
\|G_t(\mathrm{x})-G_t(\mathrm{y})\|_{\mathcal L(\R^m;X)}\le M_G\,\|\mathrm{x}-\mathrm{y}\|_X,\\
\quad
\|f_t(0)\|_X\le M_f,\quad \|G_t(0)\|_{\mathcal L(\R^m;X)}\le M_G.
\end{gathered}
\]

Fix $u\in\mathcal U$, $\mathrm{x}\in X$, and \(s\in[0,T)\). A solution of the system \eqref{ODE} with initial data $(s,\mathrm{x})$ is understood in the mild sense, namely as a continuous map $t\mapsto \Phi^u_{s,t}(\mathrm{x})$ satisfying, for each $t\in[s,T]$, the equation
\begin{equation}\label{de1}
\Phi^u_{s,t}(\mathrm{x})
=
S_{t-s}\,\mathrm{x}
+
\int_s^t S_{t-\tau}\,F_\tau\bigl(\Phi^u_{s,\tau}(\mathrm{x}),u(\tau)\bigr)\d \tau.
\end{equation}
Under hypotheses $({\bf A})$, for every triple $(s,\mathrm{x},u)\in[0,T)\times X\times\mathcal U$ there exists a unique mild solution $\Phi^u_{s,\cdot}(\mathrm{x})\in C([s,T];X)$ on the whole interval $[s,T]$; moreover, this solution is globally bounded, see \cite[Ch.~6]{Pazy1983}.

The corresponding family of maps satisfies the composition rule $\Phi^u_{t_1,t_2}\circ\Phi^u_{t_0,t_1}=\Phi^u_{t_0,t_2}$ and the normalization $\Phi^u_{t_0,t_0}=\mathrm{id}_X$ whenever $0\leq t_0\leq t_1\leq t_2\leq T$. We therefore regard $(\Phi^u_{s,t})_{0\leq s\leq t\leq T}$ as an evolution map on $X$.\footnote{Importantly, no time-reversibility is assumed.}

\subsection{Differentiability of the Mild Flow}

In addition to $({\bf A})$, we impose the following regularity assumption $({\bf A}_+)$: for a.a. $t\in I$, the maps $f_t\colon X\to X$ and $G_t\colon X\to \mathcal L(\R^m;X)$ are $C^1$. Moreover, for every $\mathrm x\in X$ the maps $t\mapsto Df_t(\mathrm x)$ and $t\mapsto DG_t(\mathrm x)$ are strongly measurable, and, for a.a. $t \in I$, 
$\mathrm x\mapsto Df_t(\mathrm x)$ and
$\mathrm x\mapsto DG_t(\mathrm x)$ are globally bounded and Lipschitz on $X$,
with a common constant $M\ge 0$.

Fix a control $\bar u\in\mathcal U$ and write
\(
\bar F_t(\mathrm x)\doteq f_t(\mathrm x)+G_t(\mathrm x)\, \bar u(t).
\)
We also abbreviate $\Delta\doteq \{(s,t)\in I^2\colon s\le t\}$.

\begin{lemma}\label{lem:diff-mild-flow}
Assume $({\bf A})$ and $({\bf A}_+)$. Then the following statements hold.

\begin{enumerate}
\item For every bounded set $B\subset X$ there exists a bounded set $\widetilde B\subset X$ such that
\(
\bar\Phi_{s,t}(\mathrm x)\in \widetilde B\)
\(\forall (s,t)\in\Delta,\ \mathrm x\in B.
\)

\item The map $(s,t,\mathrm x)\mapsto \bar\Phi_{s,t}(\mathrm x)$ is continuous on $\Delta\times X$.

\item For every $(s,t)\in\Delta$, the map $\mathrm x\mapsto \bar\Phi_{s,t}(\mathrm x)$ is Frech\'{e}t differentiable on $X$.

\item If
\(
J_{s,t}(\mathrm x)\doteq D_{\mathrm x}\bar\Phi_{s,t}(\mathrm x)\in \mathcal L(X;X),
\)
then, for every $h\in X$, the function $t\mapsto J_{s,t}(\mathrm x)h$ is the unique mild solution of the variational equation
\begin{align}
\hspace{-0.2cm} J_{s,t}(\mathrm x)\, h
&=
S_{t-s}h
+
\int_s^t
S_{t-\tau}
\Bigl(
Df_\tau(\bar\Phi_{s,\tau}(\mathrm x)) \, J_{s,\tau}(\mathrm x) h
\notag\\
&
+
DG_\tau(\bar\Phi_{s,\tau}(\mathrm x))[J_{s,\tau}(\mathrm x)\, h]\,\bar u(\tau)
\Bigr)\d \tau.
\label{variational-semigroup}
\end{align}

\item For every $h\in X$, the map
\(
(s,t,\mathrm x)\mapsto J_{s,t}(\mathrm x)\, h
\)
is continuous as $\Delta \times X \to X$.
\end{enumerate}
\end{lemma}

\subsection{Existence of Optimal Controls}\label{parag:exist}

We now turn to the existence of an optimal control. As in the ODE case, the argument has two parts: compactness of $\mathcal U$ in a weak* topology and lower semicontinuity of $\mathcal I$. A point that requires discussion is the continuity of the control-to-state map for solutions of \eqref{de1}.

We equip $\mathcal U$ with the weak* topology $\sigma(L^\infty,L^1)$ induced by the duality $L^\infty = (L^1)'$ and denote the resulting space by $\mathcal U_{w*}$. Since \(U\subset\R^m\) is compact, the set \(\mathcal U\) is bounded in \(L^\infty(I;\R^m)\) and weakly* closed. Hence, by the Banach--Alaoglu theorem, \(\mathcal U_{w*}\) is compact. Since $L^1(I;\R^m)$ is separable, this topology is metrizable on bounded sets, and hence sequential arguments are sufficient.

The second ingredient is given by the following lemma.

\begin{lemma}
Assume $({\bf A})$ and, in addition,

\begin{itemize}
\item[(${\bf B}$)] The functional \(\ell\colon X\to\R\) is lower semicontinuous.

\item[($\bf{C}$)] The operator $G$ admits the following structure:
\begin{align}
\label{G}
G_t(\mathrm x) \, \mathrm u = \sum_{j=1}^m  \bigl(\mathrm u^T \, g^j_t(\mathrm x)\bigr) \, \bm h^j,
\end{align}
where $g^j\colon I \times X \to \R^m$ satisfy the related assumptions in $({\bf A})$, and $\bm h^j\in X$. 
\end{itemize}

Let $(u^n)\subset\mathcal U$ converge to $u$ in $\mathcal U_{w*}$. Then the associated mild solutions satisfy $x^{u^n}\to x^u$ in $C(I;X)$. Moreover, the map $u\mapsto \|u\|_{L^2}^2$ is lower semicontinuous (l.s.c.) on $\mathcal U_{w*}$. Consequently, $\mathcal I$ is sequentially l.s.c. on $\mathcal U_{w*}$.
\end{lemma}

The proof is based on the mild representation \eqref{de1} and skipped for brevity: the difference $x^{u^n}-x^u$ is estimated by combining the Lipschitz bounds from $({\bf A})$ with Gr\"{o}nwall's inequality, while the control term is handled by using the finite-channel structure from $({\bf C})$ together with the weak* convergence of $(u^n)$. 

We may now apply the direct method of the calculus of variations to arrive at the following result.
\begin{theorem}
Under assumptions $({\bf A})$--$({\bf C})$, problem $(\mathbf P)$ admits a minimizer in the class $\mathcal U$.
\end{theorem}


Remark that the finite-channel structure of $G$ is essential here. Although specific, this architecture is natural for distributed systems controlled by a finite number of actuators, and still leaves enough structure to combine weak$^\ast$ compactness of controls with an exact comparison argument.

\section{Exact Cost-Increment Formula and Monotone Descent Method}

The construction developed below is local in nature and \emph{does not invoke the Hamilton--Jacobi equation}. Its central object is the family obtained by propagating the terminal cost backward along a fixed reference evolution.

Let $\bar u\in\mathcal U$ be a baseline control, $\bar\Phi\doteq\Phi^{\bar u}$ the corresponding evolution map, and $\bar x_t\doteq\bar\Phi_{0,t}(\mathrm x_0)$. Given another control $u\in\mathcal U$ and a time $s\in(0,T)$, define the element 
\begin{equation}
u\triangleright_s\bar u\doteq {\bf 1}_{[0,s)} u + {\bf 1}_{[s,T]}
\bar u \in \mathcal U \label{control-var-gen}
\end{equation}
If $\Phi\doteq\Phi^u$ denotes the evolution map generated by $u$ and $x_t\doteq x_t^u\doteq \Phi_{0,t}(\mathrm x_0)$, then
\(
x_t^{u\triangleright_s\bar u}
=
\bar\Phi_{s,t}(x_s)\),
\(0\le s\le t\le T.
\)
Now set
\(
\bar p_t \doteq \ell\circ \bar\Phi_{t,T};
\)
$\bar p_t(\mathrm x)$ is the terminal payoff obtained by starting from the state $\mathrm x$ at time $t$ with the baseline control $\bar u$. The key properties of this map are specified below. 

\begin{corollary}\label{cor:dp}
Assume $({\bf A})$, $({\bf A}_+)$, and $\ell\in C^1(X)$. 
Then, for every $t\in I$, the map $\bar p_t$ is Frech\'{e}t differentiable on $X$, and
\begin{equation}\label{Pd}
D\bar p_t(\mathrm x)
=
D\ell(\bar\Phi_{t,T}(\mathrm x))\circ D_{\mathrm x}\bar\Phi_{t,T}(\mathrm x).
\end{equation}
Moreover, for every $\eta\in X$, the map
\(
(t,\mathrm x)\mapsto D\bar p_t(\mathrm x)\,\eta
\)
is continuous on $I\times X$.
\end{corollary}

By the composition property of the evolution map, 
\begin{equation}
\bar p_t\bigl(\bar\Phi_{0,t}(\mathrm x)\bigr)=\ell\bigl(\bar\Phi_{0,T}(\mathrm x)\bigr)
\label{ppp}
\end{equation}
for every $\mathrm x\in X$ and $t\in I$. In particular, the left-hand side is independent of $t$, and therefore
\(
\ell(\bar x_T)=\bar p_T(\bar x_T)=\bar p_0(\mathrm x_0).
\)

This identity allows one to represent the terminal increment on the pair $(\bar u,u)$ in the form
\begin{equation}
\ell(x_T)-\ell(\bar x_T)
=
\bar p_T(x_T)-\bar p_0(\mathrm x_0)
=
\int_I \frac{\d}{\d t}\bar p_t(x_t)\d t.
\label{ppl}
\end{equation}
The remaining step is to compute the derivative under the integral sign. At a formal level, the variation of $t\mapsto \bar p_t(x_t)$ comes only from the mismatch between the actual control $u(t)$ and the baseline control $\bar u(t)$. Accordingly, one expects
\begin{equation}
\frac{\d}{\d t}\bar p_t(x_t)
=
(u(t)-\bar u(t))^{T}\,G_t(x_t)' \,D\bar p_t(x_t).
\label{eif}
\end{equation}
This is the key identity behind the monotone scheme.
\begin{proposition}\label{pro0}
Suppose, in addition, that $\ell$ is Lipschitz on bounded subsets of $X$. Then, the composition
$g\colon t\mapsto g(t)\doteq \bar p_t(x_t)$, $t\in I$,
is absolutely continuous on $I$, and, for a.a. $t\in I$,
\begin{equation}\label{eq:g-derivative}
g'(t)=(u(t)-\bar u(t))^T\,G_t(x_t)'D\bar p_t(x_t).
\end{equation}
\end{proposition}

\medskip
Proposition~\ref{pro0} immediately yields an exact representation for the full increment of the cost:
\[
\mathcal I[u]-\mathcal I[\bar u]
=
\int_I
\Bigl(
\bar H_t(x_t,u(t))
-
\bar H_t(x_t,\bar u(t))
\Bigr)\d t,
\]
where
\[
\bar H_t(\mathrm x,\mathrm u)
\doteq
\frac{\alpha}{2}|\mathrm u|^2
+
\mathrm u^T G_t(\mathrm x)'D\bar p_t(\mathrm x),
\ \, 
(t,\mathrm x,\mathrm u)\in I\times X\times \R^m.
\]
Thus, once the baseline control \(\bar u\) is fixed, the descent mechanism is encoded by the pointwise minimization of the reduced Hamiltonian \(\bar H_t\).

\subsection{Implementation via Sample‑and‑Hold Updates}

Take
\(
U=B_R \doteq \{\mathrm u\in\R^m\colon |\mathrm u|\le R\}.
\)
For fixed \((t,\mathrm x)\), the minimizer of
\(\mathrm u\mapsto \bar H_t(\mathrm x,\mathrm u)\) over \(U\) is
\[
\hat u_t(\mathrm x)\!=\!
\begin{cases}
\Pi_{B_R}\bigl(w_t(\mathrm x)\bigr), & \alpha>0,\\
-R\,\dfrac{G_t(\mathrm x)'D\bar p_t(\mathrm x)}
{|G_t(\mathrm x)'D\bar p_t(\mathrm x)|},
& \alpha=0,\ G_t(\mathrm x)'D\bar p_t(\mathrm x)\neq 0,
\end{cases}
\]
where \(w_t(\mathrm x)\doteq -\alpha^{-1}G_t(\mathrm x)'D\bar p_t(\mathrm x)\). 
Formally, if \(u(t)=\hat u_t(x_t)\) and \(x=x^u\), then
$$\bar H_t(x_t,u(t))\le \bar H_t(x_t,\bar u(t))\ \ \text{ for a.a. } \ \ t\in I,$$ implying that
\(\mathcal I[u]\le \mathcal I[\bar u]\).
This is a closed-loop problem. 

For $\alpha>0$, the map $(t,\mathrm x)\mapsto \hat u_t(\mathrm x)$ is continuous whenever
$(t,\mathrm x)\mapsto G_t(\mathrm x)'D\bar p_t(\mathrm x)$ is continuous. A rigorous well-posedness analysis of the resulting closed-loop equation would require additional regularity assumptions and is not pursued here.

In computations, \(D\bar p_t(\mathrm x)\) is replaced by finite differences: Let
\(
G_t(\mathrm x)\,\mathrm u
=
\sum_{j=1}^m (\mathrm u^T g_t^j(\mathrm x))\,\bm h^j.
\)
Then
\[
G_t(\mathrm x)'D\bar p_t(\mathrm x)
=
\sum_{j=1}^m D\bar p_t(\mathrm x)[\bm h^j]\,g_t^j(\mathrm x).
\]
For \(\varepsilon>0\), define
\[
\xi_j^\varepsilon(t,\mathrm x)
\doteq
\frac{
\ell\bigl(\bar\Phi_{t,T}(\mathrm x+\varepsilon \bm h^j)\bigr)
-
\ell\bigl(\bar\Phi_{t,T}(\mathrm x)\bigr)
}{\varepsilon},
\quad j=1,\dots,m.
\]
Evidently, \(\xi_j^\varepsilon(t,\mathrm x)\) approximates \(D\bar p_t(\mathrm x)[\bm h^j]\).

Fix a partition \(0=t_0<\dots<t_N=T\). Starting with \(\mathrm x^0=\mathrm x_0\), freeze the state at \(\mathrm x^k\), compute \(t\mapsto \xi^\varepsilon(t,\mathrm x^k)\) using the baseline flow \(\bar\Phi\), and set
\[
\widetilde u(t)\doteq
\Pi_{B_R}\!\left(
-\alpha^{-1}\sum_{j=1}^m \xi_j^\varepsilon(t,\mathrm x^k)\,g_t^j(\mathrm x^k)
\right),
\quad
t\in[t_k,t_{k+1}).
\]
Then solve the state equation on \([t_k,t_{k+1}]\) with
\(x_{t_k}=\mathrm x^k\) and \(u=\widetilde u\), and define
\(\mathrm x^{k+1}=x_{t_{k+1}}\).
Repeating this for \(k=0,\dots,N-1\) produces a new control on \(I\).

For sufficiently fine partitions and small \(\varepsilon\), the iteration is monotone:
\(
\mathcal I[u^{{\rm iter}+1}] \le \mathcal I[u^{\rm iter}].
\)
Indeed, as \(\max_k(t_{k+1}-t_k)\to 0\) and \(\varepsilon\to 0\), the sample-and-hold law approaches the pointwise minimizer of \(\bar H_t\), and the exact increment formula gives a nonpositive increment. In particular, the sequence \(\bigl(\mathcal I[u^{\rm iter}]\bigr)\) is nonincreasing and bounded from below by the optimal value. Hence it converges. 

\begin{algorithm}[t]
\small
\DontPrintSemicolon
\caption{Sample-and-hold update}
\label{alg:sample-hold}
Given \(u^{\rm iter}\in\mathcal U\), \(N\ge 1\), \(\varepsilon>0\)\;
\(\bar u=u^{\rm iter}\), \(\mathrm x^0=\mathrm x_0\)\;
\For{\(k=0,\dots,N-1\)}{
compute \(\xi^\varepsilon_j(\cdot,\mathrm x^k)\) from \(\bar\Phi\)\;
set \(\widetilde u(t)\!=\!\Pi_{B_R}\!(
-\alpha^{-1}\!\sum \xi_j^\varepsilon(t,\mathrm x^k)\,g_t^j(\mathrm x^k)), \ t\in [t_k,t_{k+1})\)\;
solve equation \eqref{ODE} on \([t_k,t_{k+1}]\) with \(x_{t_k}=\mathrm x^k\), \(u=\widetilde u\)\;
\(\mathrm x^{k+1}=x_{t_{k+1}}\)\;
}
\(u^{{\rm iter}+1}=\widetilde u\)\;
\end{algorithm}

\subsection{Numerical Illustration: Reaction–Diffusion Control}

Consider the semilinear reaction--diffusion equation on the one-dimensional torus $\mathbb T^1 \simeq [0,2\pi)$:\footnote{The reaction term falls outside the global Lipschitz framework of Section~II. The example is included as a numerical illustration of the descent construction in a lower-regularity setting.}
\begin{equation*}\label{eq:num-rd}
\partial_t \rho
=
\nu\,\partial_{\theta\theta}\rho
+
\beta\,\rho(1-\rho)
+
u_1(t)\,h_1(\theta)
+
u_2(t)\,h_2(\theta),
\end{equation*}
with periodic boundary conditions, where
$h_1(\theta)=\cos\theta/\sqrt{\pi}$ and
$h_2(\theta)=\sin\theta/\sqrt{\pi}$.
This model describes a population distributed along a ring-shaped habitat: diffusion models spatial dispersal and the reaction term describes local logistic growth. In the abstract notation of the semilinear system studied above, we take $x =\rho$, $X=C(\mathbb T^1)$ and define
$A=\nu\partial_{\theta\theta}$,
$f(x)=\beta x(1-x)$,
and
$G\mathrm u=\mathrm u_1 h_1+\mathrm u_2 h_2$.

The cost is the tradeoff
\begin{equation}
\mathcal I[u]
=
\frac12 \|x_T^u-\widehat x\|_{L^2(\mathbb T^1)}^2
+
\frac{\alpha}{2}\int_0^T |u(t)|^2 \,\d t.
\end{equation}
The initial and target profiles are chosen as
$x_0(\theta)=\exp(1.5\cos(\theta-1.0))$ and
$\widehat x(\theta)=\exp(2.5\cos(\theta-2.2))$.

In the reported run we use
$\nu=0.1$,
$\beta=0.05$,
$T=2$,
$\alpha=0.2$, \(R=20\),
and
$\varepsilon=10^{-3}$.
The spatial grid contains $96$ nodes, the time step is $\Delta t=10^{-3}$, the control partition has $30$ switching intervals, and the initial baseline control is $u^0\equiv 0$. We perform $4$ outer iterations of the monotone scheme. The state equation is integrated by an exponential Euler method: the linear part is treated spectrally through the Fourier representation of the heat semigroup, while the reaction and control terms are handled explicitly. 

\begin{figure}
    \centering
    \includegraphics[width=0.75\linewidth]{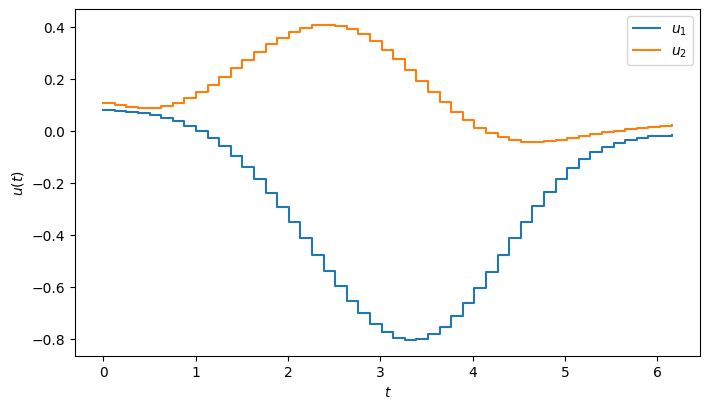}
    \includegraphics[width=0.75\linewidth]{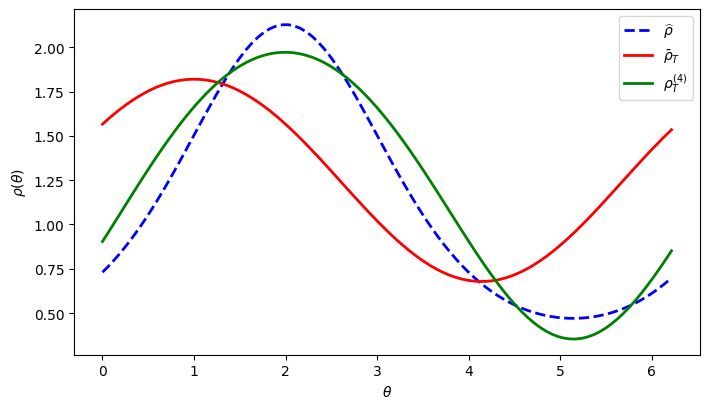}
    \caption{\small Controls generated by Alg.~1 (upper panel) and terminal distributions $\rho_T$ (lower panel) --- optimized by Alg.~1 (green) and the uncontrolled one (red) vs the target profile (dotted).  }
    \label{fig:placeholder}
\end{figure}

The minimization history of $\mathcal I$ is \(0.8283\rightarrow 0.4668\rightarrow 0.2881 \rightarrow 0.2372\rightarrow 0.2112\). Most of the decrease is achieved during the first two iterations. The terminal profile (Fig.~\ref{fig:placeholder}) moves substantially toward the target one, and the control retains the expected piecewise-constant structure.

\section{CONCLUSION}

We have obtained an exact cost-increment formula for semilinear evolution equations in Banach spaces. To the best of the authors' knowledge, this places our approach in a setting where such a construction has not been available before. Beyond the numerical direction, this formula opens the way to feedback-based variational analysis and nonclassical optimality conditions for a broader class of distributed systems.

\section*{Appendix: Proofs of Main Results}

\subsection{Proof of Lemma~\ref{lem:diff-mild-flow}}

Fix a bounded set $B\subset X$.

1, 2) The first assertion is a standard consequence of the mild flow formula, the linear growth assumption $({\bf A})$, and Gronwall's inequality. 

The second one follows from the standard comparison estimate for mild solutions corresponding to nearby triples $(s,t,\mathrm x)$, again by Gronwall's inequality.

3) Fix $(s,\mathrm x)\in[0,T)\times X$, and let $y_t\doteq \bar\Phi_{s,t}(\mathrm x)$ and $M_S\doteq \sup_{0\le r\le T}\|S_r\|_{\mathcal L(X;X)}$. 
According to Step~1, $y$ is bounded in the norm of $X$. 

For a.a. $t\in [s,T]$, we define the operator $B_t\colon X\to X$ as
\begin{align*}
v \mapsto B_t v \doteq Df_t(y_t)\, v + (DG_t(y_t) \, v)\, \bar u(t).
\end{align*}
Since $y$ is continuous on $[s,T]$, it is strongly measurable. Then, by the Carath\'{e}odory property and $({\bf A}_+)$, the maps $t\mapsto Df_t(y_t)\, v$ and $t\mapsto (DG_t(y_t)\, v)\,\bar u(t)$ are strongly measurable for any $v\in X$. Thus, $t\mapsto B_t v$ is strongly measurable. 

At the same time,
\(
\|B_t\|_{\mathcal L(X;X)}
\le
M_0 \doteq M\bigl(1+\|\bar u\|_{L^\infty}\bigr)
\)
for a.a. $t\in[s,T]$.
Hence the linear integral equation
\begin{align}
z_t
&=
S_{t-s}\eta
+
\int_s^t S_{t-\tau}B_\tau z_\tau \d\tau
\label{variational-semigroup-proof}
\end{align}
is well-posed in $C([s,T];X)$ for every $\eta\in X$. Denoting its unique solution by $z(\eta)$, one has
\(
\sup_{t\in[s,T]}\bigl\|z_t(\eta)\bigr\|_X \le M_1\|\eta\|_X
\)
(by virtue of the Gronwall inequality one may choose $M_1$ independent of $s\in[0,T]$).
Thus $\eta\mapsto z_t(\eta)$ is linear and bounded for every $t\in[s,T]$.

4) We now prove that $z_t(\eta)$ is the Fr\'echet derivative of the mild flow with respect to the initial condition.

Fix $(s,\mathrm x)\in [0,T)\times X$. Let $\eta\in X$ and set
\[
y_t^\eta \doteq \bar\Phi_{s,t}(\mathrm x+\eta),
\quad
\Delta_t^\eta \doteq y_t^\eta-y_t,
\quad
q_t^\eta \doteq \Delta_t^\eta-z_t(\eta).
\]
By item~1, for all $\eta$ in a sufficiently small ball around the origin, the trajectories $y^\eta$ remain in a common bounded subset of $X$. Moreover, the standard Lipschitz estimate for mild solutions yields
\begin{equation}
\sup_{t\in[s,T]}\|\Delta_t^\eta\|_X \le M_2\|\eta\|_X
\label{Delta-est}
\end{equation}
with a constant $M_2$ independent of $\eta$. Subtracting the mild equations for $y^\eta$ and $y$, and then subtracting \eqref{variational-semigroup-proof}, we obtain
\begin{equation}
q_t^\eta
=
\int_s^t
S_{t-\tau}
\Bigl(
B_\tau q_\tau^\eta + \rho_\tau^\eta
\Bigr)\d\tau,
\label{q-eta}
\end{equation}
\vspace{-0.2cm}
\begin{align*}
\rho_\tau^\eta
= &
R_\tau^{f,\eta}
+
R_\tau^{G,\eta}\bar u(\tau),\\
R_\tau^{f,\eta}
\doteq &
f_\tau(y_\tau^\eta)-f_\tau(y_\tau)-Df_\tau(y_\tau)\Delta_\tau^\eta,\\
R_\tau^{G,\eta}
\doteq &
G_\tau(y_\tau^\eta)-G_\tau(y_\tau)-DG_\tau(y_\tau)\Delta_\tau^\eta.
\end{align*}

For a.a. $\tau\in[s,T]$, the integral form of the mean value formula gives
\[
\begin{aligned}
&f_\tau(y_\tau^\eta)-f_\tau(y_\tau)-Df_\tau(y_\tau)\Delta_\tau^\eta
=\\
&\quad
\int_0^1
\Bigl(
Df_\tau(y_\tau+\theta\Delta_\tau^\eta)-Df_\tau(y_\tau)
\Bigr)\Delta_\tau^\eta{\d}\theta,
\end{aligned}
\]
and similarly
\[
\begin{aligned}
&G_\tau(y_\tau^\eta)-G_\tau(y_\tau)-DG_\tau(y_\tau)\Delta_\tau^\eta
=\\
&\quad
\int_0^1
\Bigl(
DG_\tau(y_\tau+\theta\Delta_\tau^\eta)-DG_\tau(y_\tau)
\Bigr)\Delta_\tau^\eta{\d}\theta.
\end{aligned}
\]
Using the Lipschitz continuity from $({\bf A}_+)$, we obtain
\[
\|R_\tau^{f,\eta}\|_X
\le
M\|\Delta_\tau^\eta\|_X^2,
\quad
\|R_\tau^{G,\eta}\|_{\mathcal L(\mathbb R^m;X)}
\le
M\|\Delta_\tau^\eta\|_X^2.
\]
Hence, by \eqref{Delta-est},
$\|\rho_\tau^\eta\|_X
\le
M M_2^2 \bigl(1+\|\bar u\|_{L^\infty}\bigr)\|\eta\|_X^2$ 
for a.a. $\tau\in[s,T]$. Therefore
\(
\|\rho^\eta\|_{L^1}
\le
M_3\|\eta\|_X^2
\)
with a constant $M_3$ independent of $\eta$.

It follows from \eqref{q-eta} that
\[
\|q_t^\eta\|_X
\le
M_S\int_s^t \|B_\tau\|_{\mathcal L(X;X)}\|q_\tau^\eta\|_X{\d}\tau
+
M_S\|\rho^\eta\|_{L^1}.
\]
Applying Gronwall's inequality, we conclude that
\[
\sup_{t\in[s,T]}\|q_t^\eta\|_X \le M_4\|\eta\|_X^2,
\]
where $M_4$ is independent of $\eta$. Thus
\[
\sup_{t\in[s,T]}
\bigl\|
\bar\Phi_{s,t}(\mathrm x+\eta)-\bar\Phi_{s,t}(\mathrm x)-z_t(\eta)
\bigr\|_X
=
o(\|\eta\|_X).
\]
Therefore $\mathrm x\mapsto \bar\Phi_{s,t}(\mathrm x)$ is Fr\'echet differentiable on $X$, and its derivative is given by
\(
D_{\mathrm x}\bar\Phi_{s,t}(\mathrm x)\,\eta = z_t(\eta).
\)
This proves item~3. Item~4 follows from the definition of $z(\eta)$.

5) Fix $\eta\in X$ and let $(s_n,t_n,\mathrm x_n)\to (s,t,\mathrm x)$ in $\Delta\times X$. 
For $r\in[s_n,T]$ set $y_r^n \doteq \bar\Phi_{s_n,r}(\mathrm x_n)$ and $z_r^n \doteq J_{s_n,r}(\mathrm x_n)\,\eta$, 
and for $r\in[s,T]$ set $y_r \doteq \bar\Phi_{s,r}(\mathrm x)$ and $z_r \doteq J_{s,r}(\mathrm x)\,\eta$.

For a.a.\ $r$ define
\begin{equation*}
\begin{aligned}
B_r^n v &\doteq Df_r(y_r^n)\, v + \big\{DG_r(y_r^n)\, v\big\}[\bar u(r)],\\
B_r v &\doteq Df_r(y_r)\, v + \big\{DG_r(y_r)\, v\big\}[\bar u(r)].
\end{aligned}
\end{equation*}

By item~4 both $z^n$ and $z$ satisfy the variational equation. Since 
$\|B_r^n\|_{\mathcal L(X;X)},\|B_r\|_{\mathcal L(X;X)}\le M_0$ for a.a.\ $r$, the estimate from item~3 gives:
\(\sup_{r\in[s_n,T]}\|z_r^n\|_X\le M_1\|\eta\|_X\) and
$\sup_{r\in[s,T]}\|z_r\|_X\le M_1\|\eta\|_X$.

If $t=s$, then $t_n-s_n\to0$ and
\begin{equation*}
z_{t_n}^n
=
S_{t_n-s_n}\eta
+
\int_{s_n}^{t_n} S_{t_n-\tau}B_\tau^n z_\tau^n{\d}\tau .
\end{equation*}
Hence $\|z_{t_n}^n-\eta\|_X\to0$, and therefore: $J_{s_n,t_n}(\mathrm x_n)\,\eta\to \eta=J_{s,s}(\mathrm x)\,\eta$.

Assume $t>s$ and set $\sigma_n\doteq\max\{s,s_n\}$. Then $\sigma_n\to s$ and, for large $n$, 
$t_n\in[\sigma_n,T]$. If $s_n\le s$, then $\sigma_n=s$ and
\begin{equation*}
z_{\sigma_n}^n
=
S_{s-s_n}\eta
+
\int_{s_n}^{s} S_{s-\tau}B_\tau^n z_\tau^n{\d}\tau .
\end{equation*}
Otherwise, $\sigma_n=s_n$ and $z_{\sigma_n}$ is defined analogously. In both cases, $\|z_{\sigma_n}^n-z_{\sigma_n}\|_X\to0$.

Define $\beta_n(\tau) = \|B_\tau^n-B_\tau\|_{\mathcal L(X;X)}$ on $[\sigma_n,T]$ and  \(\beta_n(\tau)=0\) for $\tau \in [s,\sigma_n)$. By item~2, $y_\tau^n\to y_\tau$ for a.a.\ $\tau$, hence $\beta_n(\tau)\to0$. 
Since $\beta_n(\tau)\le2M_0$, dominated convergence gives: $\|\beta_n\|_{L^1([s,T])}\to 0$. For $r\in[\sigma_n,T]$, subtract the mild flow formulas:
\begin{equation*}
\begin{aligned}
w_r^n
&\doteq z_r^n-z_r
=
S_{r-\sigma_n}(z_{\sigma_n}^n-z_{\sigma_n})\\
&+
\int_{\sigma_n}^{r} S_{r-\tau}(B_\tau^n-B_\tau)z_\tau{\d}\tau + \int_{\sigma_n}^{r} S_{r-\tau}B_\tau^n w_\tau^n{\d}\tau .
\end{aligned}
\end{equation*}
Hence
\begin{align*}
\|w_r^n\|_X
&\le
M_S\|z_{\sigma_n}^n-z_{\sigma_n}\|_X
+ M_SM_1\|\eta\|_X\!\int_s^T\!\beta_n(\tau){\d}\tau\\
&+
M_SM_0\!\int_{\sigma_n}^{r}\|w_\tau^n\|_X{\d}\tau .
\end{align*}
The first two terms tend to zero, and Gronwall's inequality yields
\(
\sup\limits_{r\in[\sigma_n,T]}\|w_r^n\|_X\to0 .
\)
In particular, $\|z_{t_n}^n-z_{t_n}\|_X\to 0$. Since $z_r$ is continuous and $t_n\to t$, also
$\|z_{t_n}-z_t\|_X\to0$. Thus
\begin{equation*}
\begin{aligned}
&\|J_{s_n,t_n}(\mathrm x_n)\,\eta-J_{s,t}(\mathrm x)\,\eta\|_X\\
&
\quad \leq \|z_{t_n}^n-z_{t_n}\|_X
+
\|z_{t_n}-z_t\|_X
\to 0.
\end{aligned}
\end{equation*}
This proves item~5.

\subsection{Proof of Corollary}

Fix $t\in I$. By item~3 of Lemma~\ref{lem:diff-mild-flow}, the map $\mathrm x\mapsto \bar\Phi_{t,T}(\mathrm x)$ is Fr\'echet differentiable on $X$ and its derivative is $J_{t,T}(\mathrm x) \doteq D_{\mathrm x}\bar\Phi_{t,T}(\mathrm x)$. Since $\ell\in C^1(X)$, the chain rule yields \eqref{Pd}. Fix $\eta\in X$; by item~5 of Lemma~\ref{lem:diff-mild-flow}, the map
\(
(t,\mathrm x)\mapsto J_{t,T}(\mathrm x)\,\eta
\)
is continuous on \(I\times X\) while $(t,\mathrm x)\mapsto \bar\Phi_{t,T}(\mathrm x)$ is continuous by item~2. Since $D\ell$ is continuous, such is the map
\(
(t,\mathrm x)\mapsto D\ell\bigl(\bar\Phi_{t,T}(\mathrm x)\bigr) J_{t,T}(\mathrm x)\, \eta
\)
 on $I\times X$. 

\subsection{Proof of Proposition~\ref{pro0}}

By Corollary~\ref{cor:dp}, $\bar p_t(\mathrm x)=\ell(\bar\Phi_{t,T}(\mathrm x))$ is $C^1$ in $\mathrm x$ for every $t\in I$, and $(t,\mathrm x)\mapsto D\bar p_t(\mathrm x)$ is continuous on $I\times X$.

\medskip

\noindent
\emph{Absolute continuity.}
For $s\in I$ set $v^s\doteq u\triangleright_s\bar u$ and $z^s\doteq x^{v^s}$. Then $g(s)=\bar p_s(x_s)=\ell(z_T^s)$.
Fix $0\le a<b\le T$. The controls $v^a$ and $v^b$ differ only on $[a,b)$.
Using the formula \eqref{de1}, the Lipschitz bounds from $({\bf A})$, boundedness of controls, and Gronwall's inequality, we get:
$\|z_T^b-z_T^a\|_X\le M_0(b-a)$
with a constant $M_0$ independent of $a,b$.
The set $\{z_T^s:s\in I\}$ is bounded, hence $\ell$ is Lipschitz on it, and
$|g(b)-g(a)|\le M_1(b-a)$
with a constant $M_1$ independent of $a,b$.
Thus $g$ is Lipschitz, hence absolutely continuous, on $I$.

\medskip

\noindent
\emph{Derivative.}
Fix $t\in[0,T)$ and $h>0$ such that $t+h\le T$. Set $y_\tau\doteq \bar\Phi_{t,\tau}(x_t)$ and $\delta_\tau\doteq x_\tau-y_\tau$ for $\tau\in[t,t+h]$. Then $\delta_t=0$.
By the composition property of $\bar\Phi$,
\begin{equation*}
\bar p_t(x_t)
=\ell(\bar\Phi_{t,T}(x_t))
=\ell(\bar\Phi_{t+h,T}(y_{t+h}))
=\bar p_{t+h}(y_{t+h}),
\end{equation*}
hence
\begin{equation}\label{eq:g-diff-1}
g(t+h)-g(t)=\bar p_{t+h}(x_{t+h})-\bar p_{t+h}(y_{t+h}).
\end{equation}
Since $\bar p_{t+h}$ is Fr\'echet differentiable at $y_{t+h}$, we have
\begin{equation*}\label{eq:g-diff-2}
\begin{gathered}
\bar p_{t+h}(x_{t+h})-\bar p_{t+h}(y_{t+h})
=D\bar p_{t+h}(y_{t+h})[\delta_{t+h}]+r_h,
\end{gathered}
\end{equation*}
where \(r_h=o(\|\delta_{t+h}\|_X)\). Subtracting the equations for $x$ and $y$ on $[t,t+h]$ gives
\begin{align}
&\delta_{t+h}
=\int_t^{t+h} S_{t+h-\tau}\Bigl(
f_\tau(x_\tau)-f_\tau(y_\tau)\label{eq:delta-mild}
\\
& \ \ +\bigl(G_\tau(x_\tau)-G_\tau(y_\tau)\bigr) \, u(\tau)
+G_\tau(y_\tau)\, (u(\tau)-\bar u(\tau))
\Bigr)\d\tau.\notag
\end{align}
By $({\bf A})$ and boundedness of $u$, Gronwall's inequality yields
\begin{equation}\label{eq:delta-bound}
\sup_{\tau\in[t,t+h]}\|\delta_\tau\|_X\le M_t h,
\end{equation}
hence $\|\delta_{t+h}\|_X=O_t(h)$ and therefore $r_h=o_t(h)$.
Using the continuity of $(s,\mathrm x)\mapsto D\bar p_s(\mathrm x)$ and $y_{t+h}\to x_t$, we also have
\begin{equation}\label{eq:dpt-shift}
D\bar p_{t+h}(y_{t+h})[\delta_{t+h}]
=D\bar p_t(x_t)[\delta_{t+h}]+o_t(h).
\end{equation}
Combining \eqref{eq:g-diff-1}--\eqref{eq:dpt-shift} yields
\begin{equation}\label{eq:g-diff-3}
g(t+h)-g(t)=D\bar p_t(x_t)[\delta_{t+h}]+o_t(h).
\end{equation}

Now decompose \eqref{eq:delta-mild} as
\begin{equation}\label{eq:delta-split}
\delta_{t+h}
=\int_t^{t+h}S_{t+h-\tau}G_t(x_t)\,\bigl(u(\tau)-\bar u(\tau)\bigr) \d\tau+\rho_h,
\end{equation}
where $\rho_h$ is the sum of the remaining terms, namely the integrals with
$f_\tau(x_\tau)-f_\tau(y_\tau)$,
$\bigl(G_\tau(x_\tau)-G_\tau(y_\tau)\bigr)u(\tau)$,
and $\bigl(G_\tau(y_\tau)-G_t(x_t)\bigr)(u(\tau)-\bar u(\tau))$.

Using \eqref{eq:delta-bound} and the Lipschitz bounds from $({\bf A})$, the first two parts give $\|\cdot\|_X$-contributions of order $O_t(h^2)$. For the last part, split: \[
\begin{aligned}
&\|G_\tau(y_\tau)-G_t(x_t)\|\\
&\quad \le
\|G_\tau(y_\tau)-G_\tau(x_t)\|+\|G_\tau(x_t)-G_t(x_t)\|.
\end{aligned}
\]
The first term contributes $o_t(h)$ because $y_\tau\to x_t$ as $\tau\downarrow t$ and
\begin{equation*}
\int_t^{t+h}\|G_\tau(y_\tau)-G_\tau(x_t)\|\d\tau
\le
M_G\int_t^{t+h}\!\!\!\|y_\tau-x_t\|_X\d\tau. 
\end{equation*}

To handle the second term, note that the set $K\doteq\{x_t:t\in I\}$ is compact in the strong topology of $X$.
Choose a dense subset $(\mathrm x^k)_{k\ge1}\subset K$.
For each $k$, the map $\tau\mapsto G_\tau(\mathrm x^k)\in\mathcal L(\mathbb R^m;X)$ is strongly measurable by $({\bf A})$,
hence by the Lebesgue differentiation theorem there exists a set $J_k\subset I$ of full measure such that
\begin{equation*}
\frac{1}{h}\int_t^{t+h}\bigl\|G_\tau(\mathrm x^k)-G_t(\mathrm x^k)\bigr\|\d\tau\to 0
\ \text{ as }h\downarrow0
\end{equation*}
for every $t\in J_k$.
Let $J\doteq\bigcap_{k\ge1}J_k$, then $J$ has full measure.

Fix $t\in J$. For any $\varepsilon>0$ pick $k$ with $\|x_t-\mathrm x^k\|_X\le\varepsilon$.
By the Lipschitz estimate of $G_\tau$ in $({\bf A})$, for $h>0$,  we get
\begin{align*}
&\frac{1}{h}\int_t^{t+h}\|G_\tau(x_t)-G_t(x_t)\|\d\tau
\\[-0.13cm]
&\quad\le
\frac{1}{h}\int_t^{t+h}\|G_\tau(x_t)-G_\tau(\mathrm x^k)\|\d\tau\\[-0.13cm]
&\quad+
\frac{1}{h}\int_t^{t+h}\|G_\tau(\mathrm x^k)-G_t(\mathrm x^k)\|\d\tau\\[-0.13cm]
&\quad +
\frac{1}{h}\int_t^{t+h}\|G_t(\mathrm x^k)-G_t(x_t)\|\d\tau\\[-0.cm]
&\le
2M_G\varepsilon
+
\frac{1}{h}\int_t^{t+h}\|G_\tau(\mathrm x^k)-G_t(\mathrm x^k)\|\d\tau.
\end{align*}
Letting $h\downarrow0$ and then $\varepsilon\downarrow0$ yields
\begin{equation*}
\int_t^{t+h}\|G_\tau(x_t)-G_t(x_t)\|\d\tau=o_t(h)
\quad\forall \, t\in J.
\end{equation*}
Therefore $\|\rho_h\|_X=o_t(h)$ for every $t\in J$, hence also
\begin{equation}\label{eq:dpt-rho}
D\bar p_t(x_t)[\rho_h]=o_t(h).
\end{equation}

Substituting \eqref{eq:delta-split} into \eqref{eq:g-diff-3} and using \eqref{eq:dpt-rho}, we obtain
\begin{equation}\label{eq:g-diff-4}
\begin{aligned}
&\frac{g(t+h)-g(t)}{h}
=o_t(1) \ +\\
&\quad\frac{1}{h}\int_t^{t+h} \!\!\!
D\bar p_t(x_t)\bigl[S_{t+h-\tau}G_t(x_t)\, (u(\tau)-\bar u(\tau))\bigr]\d\tau.
\end{aligned}
\end{equation}
Set $w_t\doteq G_t(x_t)'D\bar p_t(x_t)\in\mathbb R^m$. Then
$$D\bar p_t(x_t)[G_t(x_t)v]=v^T w_t$$ for all $v\in\mathbb R^m.$
Since $\mathbb R^m$ is finite-dimensional and $(S_\sigma)_{\sigma\ge0}$ is strongly continuous, the convergence $S_\sigma\to \mathrm{id}_X$ is uniform on the compact set $G_t(x_t)\{v:|v|\le R\}$ for any fixed $R>0$.
As the controls are bounded, $|u(\tau)-\bar u(\tau)|$ is uniformly bounded, hence in \eqref{eq:g-diff-4} one may replace $S_{t+h-\tau}$ by $\mathrm{id}_X$ at the price of $o_t(1)$ after division by $h$. Therefore
\begin{equation*}
\frac{g(t+h)-g(t)}{h}
=
\frac{1}{h}\int_t^{t+h}(u(\tau)-\bar u(\tau))^T w_t\d\tau
+o_t(1).
\end{equation*}
If $t$ is also a Lebesgue point of $u$ and $\bar u$, then
\begin{equation*}
\frac{1}{h}\int_t^{t+h}(u(\tau)-\bar u(\tau))\d\tau\to u(t)-\bar u(t)
\ \text{ as }h\downarrow 0,
\end{equation*}
and hence \eqref{eq:g-derivative} follows for such $t$.
Since $g$ is a.e.\ differentiable, the Lebesgue point property for $u,\bar u$ holds a.e., and $J$ has full measure, \eqref{eq:g-derivative} holds for a.a.\ $t\in I$.

\bibliographystyle{IEEEtran}
\bibliography{starmax_en_cleaned_sorted}


\end{document}